\documentclass[10pt]{article}
\usepackage[utf8]{inputenc}
\usepackage{tabularx,lscape,longtable}
\usepackage[english]{babel}
\usepackage[dvipsnames]{xcolor}
\usepackage{amsfonts,comment}
\usepackage{multicol}
\usepackage{amssymb,latexsym,amsmath}
\makeatletter \oddsidemargin 0in \evensidemargin 0in \textwidth
16cm \RequirePackage[dvips]{graphicx} \textheight 20cm
\DeclareGraphicsExtensions{.jpg, .eps}
\newtheorem{t1}{Theorem}[section]
\newtheorem{p1}{Proposition}[section]
\newtheorem{l1}{Lemma}[section]
\newtheorem{c1}{Corollary}[section]
\newtheorem{d1}{Definition}[section]
\newtheorem{r1}{Remark}[section]

\newtheorem{ex}{Example}[section]
\definecolor{cadmiumred}{rgb}{0.89,0.0,0.13}
\renewcommand{\baselinestretch}{1.3}
\renewcommand{\theequation}{\thesection.\arabic{equation}}
\begin{document}
\title{\bf
Use of copula functions  in  error assessment due to deviation from dependence assumption}
{\author{Subarna Bhattacharjee$^{1}$$^{}$\thanks{Corresponding author~:
E-mail ID: subarna.bhatt@gmail.com}, Aninda Kumar Nanda$^{2},$ Subhashree Patra$^{3}$\\
{\it $^{1,3}$ Department of Mathematics, Ravenshaw University, Cuttack-753003, Odisha, India}\\
{\it $^{2}$ Indian Statistical Institute, Delhi Centre, India.}\\}
\date{\today}
\maketitle
\begin{abstract}
In this paper, we analyze the relative errors in various reliability measures due to the tacit assumption that the components associated with a $n$-component series system or a parallel system are independently working where the components are dependent. We use  Copula functions in said error analysis. This technique generalizes the existing work on error assessment for many wide class of distributions.
\end{abstract}
\section{Introduction}
The analysis of error incurred in computation of various reliability functions arising due to
an erroneous assumption of components working independently has been emphasized
by various authors in recent time. The articles by Klein and Moeschberger (1986, 1987),
Moeschberger and Klein (1984), Gupta and Gupta (1990), and Nanda et al.(2022) focus on
the relative errors of some bivariate exponential distributions.
In this paper, we discuss the series of $n$- components and the parallel system and do the error analysis using copula functions. The previous work in this direction as listed above has been done for specific bivariate distributions. However, we investigate the errors in reliability functions using copula functions since every joint distribution can be expressed in terms of some specific copula. This will help us to generalize the results obtained by aforementioned authors. We establish some results related to error assessment, particularly over- and under-evaluation in various reliability functions like survival functions (SF), mean residual life function (MRL), hazard rate (HR), reversed hazard rate (RHR), aging intensity function (AI), for broader class of distributions. As usual, OA (UA) of a aging function refers to the difference of said aging function with independent from dependent counterpart being negative (positive). \\

We recall that, for a non-negative random variable $T$ having density function $f(\cdot),$ and survival function (SF) $\bar{F}(\cdot),$ the aging notions  are described at an arbitrary time $t$ through hazard rate (HR)  $r(t)=f(t)/\bar{F}(t),$ reversed hazard rate (RHR) $\mu(t)=f(t)/{F}(t),$ mean residual life function (MRL) $m(t)=E(T-t \mid T>t)=\frac{\int_{t}^{\infty} \bar{F}(x)dx}{\bar{F}(t)}$ and aging intensity function (AI) $L(t)=tr(t)/\int_{t}^{\infty} r(u)du.$\\
${}$\hspace{0.8cm}Here, we give a brief overview of copula function. According to Sklar (1959), for any $n$ dimensional joint cumulative (survival) function $F$ ($\bar{F}$) with marginal distribution (survival) functions $F_i$ ($\bar{F}_{i}$), there exists a function, say $C$ ($\hat{C}$), called copula (survival copula) which couples the joint distribution (survival) function with their marginal distributions (survival) functions, i.e., for a copula $C(\cdot)$, we have $$F_{T_{1}, T_{2}, \ldots, T_{n}}(t_1,t_2,\ldots,t_n)=P(T_{1}<t_1,T_2<t_2,\ldots,T_n<t_n)=C(F_{1}(t_1), F_{2}(t_2), \ldots, F_{n}(t_n)),$$ and for survival copula we have $$\bar{F}_{T_{1}, T_{2}, \ldots, T_{n}}(t_1,t_2,\ldots,t_n)=P(T_{1}>t_1,T_2>t_2,\ldots,T_n>t_n)=\hat{C}(\bar{F}_{1}(t_1), \bar{F}_{2}(t_2), \ldots, \bar{F}_{n}(t_n)),$$ for $t\geq 0.$\\
${}$\hspace{0.8cm}In section 2, we give the expressions of survival function of $n$ component series and parallel systems in terms of copula (survival copula) function. We establish that, for a given distribution, we need to arrange $\bar{F}_{S}^{I}(t)$ and  $\bar{F}_{S}^{D}(t)$ in increasing order at an arbitrary $t,$ and so is for $\bar{F}_{P}^{I}(t)$ and  $\bar{F}_{P}^{D}(t).$ Further, we show that over assessment and under assessment of different reliability functions can be inferred on the basis of copula functions. The magnitude of error or relative error can be expressed in terms of copula functions. The role of stochastic orders in knowing the error assessment is discussed in Table \ref{tab1} and we demonstrate it for some well-known copula functions. The findings of this article are discussed in Section 3. 
\section{Copula functions in error analysis}
We note that the survival function, denoted by $\overline{F}_{S}(\cdot)$ of the lifetime $T_S (=\min_{1\leq i\leq n} T_i)$ of a series system formed by $n$ components with $ith$ component having lifetime $T_{i}$ having distribution function $F_{i}(\cdot)$ for $i=1,2,\ldots,n$ is given by \begin{equation}
\label{eq1s}
\overline{F}_{S}(t)
=\overline{F}(t,t,\ldots,t)=\hat{C}(\overline{F}_{1}(t),\overline{F}_{2}(t),\ldots,\overline{F}_{n}(t))=\hat{C}(\hat{ u}(t))
\end{equation}
where $\hat{u}(t)$ is the vector with components consisting of marginal distribution functions, $\hat{u}_i(t)=\bar{F}_{i}(t),$ i.e., $$\hat{u}(t)=(\hat{u}_1(t),\hat{u}_2(t),\ldots,\hat{u}_n(t)), $$ here, $\bar{F}$ represents the joint survival function of $T_1,T_2,\ldots,T_n,$ and $\hat{C}$ is the corresponding survival copula. Similarly,
the survival function, denoted by $\bar{F_P}(\cdot)$  of an $n$ component parallel system with lifetime $T_P (=\max_{1\leq i\leq n} T_i)$ is given by
\begin{equation}
\label{eq2p}
\overline{F}_{P}(t)=1-F(t,t,t,\cdots,t)=1-C(F_{1}(t), F_{2}(t), \ldots, F_{n}(t))=1-C(u(t)),
\end{equation}
where $C$ is the copula function, and the vector consisting of marginal distribution function is ${u}(t)=({u}_1(t),{u}_2(t),\ldots,{u}_n(t)),$ with ${u}_i(t)={F}_{i}(t), t\geq 0.$
In other words, the distribution of $T_P$, representing the lifetime of a parallel system having $n$ components is given by
$F_P(t)={C}({F}_{1}(t),{F}_{2}(t),\ldots,{F}_{n}(t))={C}(u).$ 
 \\

${}$\hspace{0.5cm}In subsequent discussions, we denote the lifetime of a series system having dependent and independent components by $T_{S}^{I},$ and $T_{S}^{D}$ respectively. For  parallel system, we denote them by $T_{P}^{I},$ and $T_{P}^{D}$ respectively. The survival function of series and parallel system are denoted by $\overline{F}_{S}^{I}(\cdot)$ and $\overline{F}_{P}^{I}(\cdot)$ respectively if its components are independently working. For systems with dependent components, we denote the corresponding functions by $\overline{F}_{S}^{D}(\cdot)$ and $\overline{F}_{P}^{D}(\cdot)$ respectively. We use analogous notations for other aging functions namely hazard rate, and reversed hazard rate.
In case, the $n$ components in a parallel system work independently, then $F^{I}_P(t)={C}(u)=\prod_{i=1}^{n}u_{i}(t),$ whereas for a series system we have $\bar{F}^{I}_S(t)={\hat{C}}(u)=\prod_{i=1}^{n}\hat{u}_{i}(t),$ $t\geq 0.$\\

${}$\hspace{0.5cm}Nadarajah et al. (2017) has given a detailed list of parametric copulas. In existing literature on said error assessment, we find that different authors have studied the OA and UA in error analysis of reliability functions for specific distributions say Gumbel-I, Gumbel-II, Gumbel-III, Freund, Sarkar etc. In this paper, we establish that the said analysis can be alternatively done using copula (survival copula) functions, which thereby generalize results for wider class of parametric distributions. In Table \ref{tab1}, we list out a few copulas and find the corresponding relationships between dependent and independent counterparts of series and parallel systems. 
\subsection{Relation between survival function of  systems having dependent or independent components}
${}$\hspace{0.5cm}For larger interest of reliability practitioners, we state that the survival function of any $n$ component series system at an arbitrary time point is less than the survival function of any $n$ component parallel system at the specified time point irrespective of the fact that the components are independent or dependent. But, in general, we cannot order $\overline{F}_{S}^{I}(\cdot)$ and $\overline{F}_{S}^{D}(\cdot)$ for an arbitrary multi-variate distribution. Similar situations arises when we consider $\overline{F}_{P}^{I}(\cdot)$ and $\overline{F}_{P}^{D}(\cdot).$ The error assessment of a series (parallel) system corresponding to a particular function is all about analyzing the difference between dependent and independent counterpart.\\

${}$\hspace{0.5cm}In upcoming theorem, we prove that for any  multivariate probability distribution, there is stochastic ordering between the lifetime of $n$ component series and parallel systems having dependent or independent components. We infer that $T_{S}$ and $T_{P}$ are ordered in usual stochastic sense under two different conditions, namely, when components are dependently or independently working. To the best of my knowledge, no such result exists in literature but it may seem to be very obvious for many readers. An innovative approach of proving the following theorem is presented.
\begin{t1}
\label{multi}
For any multivariate distribution,
\begin{enumerate}
\item[(i)] $\overline{F}_{P}^{I}(t) \geq \overline{F}_{S}^{I}(t),$ i.e., $1-\prod_{i=1}^{n}P(T_i<t) \geq \prod_{i=1}^{n} P(T_i>t);$
\item[(ii)]  $\overline{F}_{P}^{I}(t) \geq \overline{F}_{S}^{D}(t),$ i.e., $1-\prod_{i=1}^{n}P(T_i<t)\geq P(T_1>t,T_2>t, \ldots,T_n>t);$
\item[(iii)]  $\overline{F}_{P}^{D}(t) \geq \overline{F}_{S}^{I}(t),$ i.e., $1-P(T_1<t,T_2<t, \ldots,T_n<t)\geq \prod_{i=1}^{n}P(T_i>t);$
\item[(iv)] $\overline{F}_{P}^{D}(t) \geq \overline{F}_{S}^{D}(t),$ i.e., $1-P(T_1<t,T_2<t, \ldots,T_n<t)\geq P(T_1>t,T_2>t, \ldots,T_n>t) ;$
\end{enumerate}
\end{t1}
{\bf Proof.} For some $k\in \{1,2,\ldots,n\},$ we have $$\overline{F}_{P}^{I}(t)=1-\prod_{i=1}^{n}F_{i}(t)\geq 1-F_{k}(t)=\overline{F}_{k}(t) \geq \prod_{i=1}^{n} \overline{F}_{i}(t)=\overline{F}_{S}^{I}(t),$$
 proving $(i).$ However, \begin{eqnarray}
\overline{F}_{P}^{I}(t)-\overline{F}_{S}^{D}(t)&=&1-\Big\{\prod_{i=1}^{n} F_{i}(t)+
\overline{F}(t,t,\ldots,t)\Big\}\nonumber\\
&=&1-\Big\{\prod_{i=1}^{n}P(A^{c}_{i})+P\big(\bigcap_{i=1}^{n}A_{i}\big)\Big\}\nonumber\\&=&1-\Big\{\prod_{i=1}^{n}P(A^{c}_{i})+P\Big(\big(\bigcup_{i=1}^{n}A_{i}^{c}\big)^{c}\Big)\Big\}\nonumber\\
&=&P\big(\bigcup_{i=1}^{n} A_{i}^{c}\big)-\prod_{i=1}^{n}P\big(A_{i}^{c}\big)\nonumber\\
&=&P\big(\bigcup_{i=1}^{n} A_{i}^{c}\big)-P\big(A_{i}^{c}\big)
\geq 0,\nonumber
\end{eqnarray}
proving $(ii).$ Proof of $(iii)$ can be done on a similar line. Also, $\overline{F}_{P}^{D}(t)=1-F(t,t,\ldots)\geq 1-F_{k}(t)=\overline{F}_{k}(t) \geq \prod_{i=1}^{n}\overline{F}_{i}(t)$ For all $i=1,2,\ldots,n,$ defining events $A_{i}$ as $A_{i}=(X_i>t)$ we note that
\begin{eqnarray}
\overline{F}_{P}^{D}(t)-\overline{F}_{S}^{D}(t)&=&1-\Big\{
F(t,t,\ldots,t)+\overline{F}(t,t,\ldots,t)\Big\}\nonumber\\
&=&1-\Big\{P(\bigcap_{i=1}^{n}A^{c}_{i})+P\big(\bigcap_{i=1}^{n}A_{i}\big)\Big\}\nonumber\\
&=&1-\Big\{P\big(\big(\bigcup_{i=1}^{n} A_{i}\big)^{c}\big)+P\big(\bigcap_{i=1}^{n}A_{i}\big)\Big\}\nonumber\\
&=&P\big(\bigcup_{i=1}^{n} A_{i}\big)-P\big(\bigcap_{i=1}^{n}A_{i}\big)
\geq 0,\nonumber
\end{eqnarray}
the last inequality follows since $\bigcap_{i=1}^{n}A_{i} \subset \bigcup_{i=1}^{n}A_{i}.$ This proves $(ii).$ $\hfill\Box$\\
The following remark gives an alternative way to prove $(i)$ of Proposition \ref{multi}.
\begin{r1}
If $T_{i}'$ s are independent random variables then $T_{S}\leq_{hr} T_{P},$ i.e., $r_{S}^{I}(t)\geq r_{P}^{I}(t)$ for all $t>0.$
\end{r1}
{\bf Proof.} If $T_{i}'$ s are independent random variables then $T_{(k)}\leq_{hr}T_{(k+1)},$ for $k\in\{1,2,\ldots,n-1\},$  where $T_{(i)}$ is the $i$th order statistic (cf. Shaked and Shanthikumar (2007) ). This gives $(T_{P}=)T_{(1)}\leq_{hr}T_{(n)}(=T_{S}),$ provided the components are independently working. This completes the proof.$\hfill\Box$
\begin{r1}
Readers may note that, in Theorem \ref{multi}, $(i)$ is a special case of $(iii),$ and $(i),$ $(ii)$ and $(iii)$ are particular cases of $(iv)$.  In a nutshell, from $(iv)$ one can arrive at other results. But, one may look at the independent proofs of $(i), (ii)$ and $(iii).$
\end{r1}
\begin{r1}
The two pair of functions satisfy $(\bar{F}_{S}^{I}(t),\bar{F}_{S}^{D}(t))\leq (\bar{F}_{P}^{I}(t),\bar{F}_{P}^{D}(t)),$ where the notation $(x,y) \leq (z,w)$ means $x\leq z, x\leq w, y\leq z$ and $y\leq w.$ So, for a distribution, we need to find out how $\bar{F}_{S}^{I}(t)$ and $\bar{F}_{S}^{D}(t)$ are ordered or note whether they cross each other. The same logic holds for the parallel system.
\end{r1}
\subsection{Expressions relating multivariate survival copula, copula and incurred errors}
An independent approach to establish the relation between joint survival and joint distribution of $n$ random variables is derived in the following proposition, which also relates a copula and the corresponding survival copula. Readers are quite familiar with bivariate result, namely, $$\hat{C}(\hat{u}_{1},\hat{u}_{2})=1-u_1-u_2+C(u_1,u_2).$$ which follows as a particular case from upcoming proposition. In the literature, we find mostly examples of copulas. Moreover, the following proposition will also help us in answering to a  question that often may ponder in the minds of readers as to how to obtain survival copula for a given copula and vice-versa.  Interestingly, for radially symmetric copulas, for which $C\equiv \hat{C},$ replacing $u (\hat{u})$ by $(1-u) (u)$ in a copula $C(u) (\hat{C}(\hat{u}))$, one may get $\hat{C}(\hat{u}) (C(u))$.
\begin{p1}
The relation between joint survival and joint distribution function of  random variables $X_1,X_2,\ldots,X_n,$ is given by
$$P(X_1>t_1,X_2>t_2, \ldots, X_n>t_n)=1-S_1+S_2-S_3+\ldots-(-1)^{n-1}P(X_1<t_1,X_2<t_2,\ldots,X_n<t_n)$$
where $$S_1=\sum_{i=1}^{n}P(X_i<t_i), S_2=\mathop{\sum_{i=1}^{n}}_{i\neq j}\sum_{j=1}^{n}P(X_i<t_i,X_j<t_j ),\ldots,S_n=P(X_1<t_1,X_2<t_2,\ldots,X_n<t_n).$$
In other words, $$\bar{F}_{X_1X_2, \ldots X_n}(t_1,t_2, \ldots, t_n)=1-S_1+S_2-S_3+\ldots-(-1)^{n-1}F_{X_1X_2 \ldots X_n}(t_1,t_2, \ldots, t_n)$$ i.e., for $\hat{u}=(\hat{u}_1,\hat{u}_2,\ldots,\hat{u}_n)$ and $u=(u_1,u_2,\ldots,u_n),$ such that $\hat{u}_i=1-u_i,$
we have
$$\hat{C}(\hat{u})=1-S_1+S_2-S_3+\ldots-(-1)^{n-1}S_{n}.$$
\end{p1}
{\bf Proof.} We use Poincare' theorem to prove the theorem which states that if $A_1,\ldots,A_n$ are events not necessarily mutually exclusive then the probability that at least one of $A_1, \ldots,A_n$ occurs is
\begin{equation}
\label{poncare}
   P(\cup_{i=1}^{n}A_i)=S_1-S_2+S_3-\ldots+(-1)^{n-1}S_n, 
\end{equation} where 
$S_1=\sum_{i=1}^{n}P(A_i),$ $S_2=\mathop{\sum_{i=1}^{n}}_{i\neq j}\sum_{j=1}^{n}P(A_iA_j),\ldots, S_n=P(A_1A_2\ldots A_n).$\\

Considering the events $A_i$ equal to $(X_i<t_i)$ for $i=1,2,\ldots,n$
and using (\ref{poncare}), we get
$$P(\cap A_{i}^{c})=P(\cup A_i)^{c} =1-P(\cup A_i)=1-S_1+S_2-S_3+\ldots-(-1)^{n-1}S_n,$$ i.e.,
$$P(X_1>t_1,X_2>t_2, \ldots, X_n>t_n)=1-S_1+S_2-S_3+\ldots-(-1)^{n-1}P(X_1<t_1,X_2<t_2,\ldots,X_n<t_n).$$ This proves the theorem. $\hfill\Box$\\

We know that the relative error incurred in survival function of any system due to the assumption of independence among components is the ratio of the difference of survival function having independent components from survival function having dependent components and the survival function having independent components of the corresponding system. Similarly, the relative error is defined for other aging functions, taking into account the corresponding counterparts. Thereby, the error, denoted by $E^{\overline{F}}_{S}(t)$, incurred in the survival function due to the incorrect assumption that the dependent components $n$ forming a series system are independently working is expressed in terms of the survival copula as \begin{eqnarray}
 \label{series}
E^{\overline{F}}_{S}(t)&=&\frac{\bar{F}^{D}_{S}(t)-\bar{F}^{I}_{S}(t)}{\bar{F}^{I}_{S}(t)}\nonumber\\
&=&\frac{\overline{F}(t,t,\ldots,t)-\prod_{i=1}^{n}\overline{F}_{i}(t)}{\prod_{i=1}^{n}\overline{F}_{i}(t)}\nonumber\\
&=&\frac{\hat{C}(u)-\prod_{i=1}^{n}\hat{u}_i(t)}{\prod_{i=1}^{n}\hat{u}_i(t)}\nonumber
\end{eqnarray}

Similarly, the error, denoted by $E^{\overline{F}}_{P}(t)$, incurred in survival function due to the wrong assumption that $n$ dependent components forming a parallel system are independently working is expressed in term of copula as \begin{eqnarray}
    \label{parallel}E^{\overline{F}}_{P}(t)&=&\frac{\bar{F}^{D}_{P}(t)-\bar{F}^{I}_{P}(t)}{\bar{F}^{I}_{P}(t)}\nonumber\\
&=&\frac{\prod_{i=1}^{n}{F}_{i}(t)-{F}(t,t,\ldots,t)}{1-\prod_{i=1}^{n}{F}_{i}(t)}\nonumber\\
&=&\frac{\prod_{i=1}^{n}u_{i}(t)-C(u)}{1-\prod_{i=1}^{n}u_{i}(t)}.\nonumber
\end{eqnarray}
The over-assessment (OA) or under assessment (UA) in a reliability function, say survival function, hazard rate, reversed hazard rate arising due to the wrong assumption that components are independently working when they are actually dependent is judged on the basis of sign of the difference of concerned function of dependent and independent counterpart appearing in the numerator. Clearly, for OA (UA) of a reliability function, the relative error (or error) is negative (positive).  It suffices to focus on the sign of error in a  function rather than on its relative error to infer about OA or UA of the corresponding function. \\
For the sake of simplicity, we focus on the error of failure rate function
for a $n$ component series system (i.e. numerator of corresponding relative error) and obtain
\begin{eqnarray}
r_{S}^{D}(t)-r_{S}^{I}(t)&=& -\frac{d}{dt} \ln \overline{F}(t,t,\ldots,t)-\sum_{i=1}^{n}r_{T_{i}}(t)\nonumber\\
   &=& -\frac{d}{dt} \ln \hat{C}\Big(\overline{F}_{T_1}(t),\overline{F}_{T_2}(t),\ldots,\overline{F}_{T_n}(t))\Big)+\sum_{i=1}^{n}\frac{d}{dt}\ln \overline{F}_{T_i}(t)\nonumber\\
&=& -\frac{d}{dt} \Big\{\ln \hat{C}\Big(\overline{F}_{T_1}(t),\overline{F}_{T_2}(t),\ldots,\overline{F}_{T_n}(t))\Big)-\sum_{i=1}^{n}\ln \overline{F}_{T_i}(t)\Big\}\nonumber\\
&=& -\frac{d}{dt} \ln \Big\{\frac{\hat{C}\Big(\overline{F}_{T_1}(t),\overline{F}_{T_2}(t),\ldots,\overline{F}_{T_n}(t))\Big)}{\prod_{i=1}^{n}\overline{F}_{T_{i}}(t)}\Big\}\nonumber\\&=&-\frac{d}{dt} \ln \frac{\overline{F}(t,t,\ldots,t)}{\prod_{i=1}^{n}\overline{F}_{T_{i}}(t)}\nonumber\\&=&-\frac{d}{dt}\ln \Big(\frac{\hat{C}(\hat{u})}{\prod_{i=1}^{n}\hat{u}_{i}(t)}\Big)\nonumber
\end{eqnarray}
where $r_{S}^{D}(t)$ and $r_{S}^{I}(t)$ represent the failure rate of $n$ component series system having dependent and independent components respectively.
For the reversed hazard rate of a series system, we have
\begin{eqnarray}\mu_{S}^{D}(t)-\mu_{S}^{I}(t)
&=&\frac{d}{dt}\ln \Big(1-\overline{F}(t,t,\ldots,t)\Big)-\frac{d}{dt}\ln \Big(1-\prod_{i=1}^{n}\overline{F}_{T_i}(t)\Big)\nonumber\\
&=&\frac{d}{dt}\ln \Big(\frac{1-\overline{F}(t,t,\ldots,t)}{1-\prod_{i=1}^{n}\overline{F}_{T_i}(t)}\Big)\nonumber\\&=&\frac{d}{dt} \ln \Big(\frac{1-\hat{C}(\hat{u})}{1-\prod_{i=1}^{n}\hat{u}_{i}(t)}\Big),\nonumber
\end{eqnarray}
where $\mu_{S}^{D}(t)$ and $\mu_{S}^{I}(t)$ represent the reversed hazard rate of $n$ component series system having dependent and independent components respectively.\\

For parallel system, using notations of similar form, we obtain
\begin{eqnarray}
\overline{F}_{P}^{D}(t)-\overline{F}^{I}_{P}(t)&=&\big(1-F(t,t,\ldots,t)\big)-\big(1-\prod_{i=1}^{n}F_{T_{i}}(t)\big)\nonumber\\
&=&\Big(\prod_{i=1}^{n}u_{i}(t)\Big)-C(u), \nonumber
\end{eqnarray}
giving
\begin{eqnarray}
r^{D}_{P}(t)-r^{I}_{P}(t)&=&-\frac{d}{dt} \ln \Big(\frac{1-F(t,t,\ldots,t)}{1-\prod F_{T_i}(t)}\Big)\nonumber\\
&=&-\frac{d}{dt}\ln \Big(\frac{1-C(u)}{1-\prod_{i=1}^{n}u_{i}(t)}\Big),\nonumber
\end{eqnarray}
and
\begin{eqnarray}
\mu^{D}_{P}(t)-\mu^{I}_{P}(t)&=&\frac{d}{dt}\ln \Big(\frac{F(t,t,\ldots,t)}{\prod F_{T_i}(t)}\Big)\nonumber\\&=&\frac{d}{dt}\ln\Big(\frac{C(u)}{\prod u_{i}(t)}\Big).\nonumber
\end{eqnarray}

 \subsection{Establishing OA and UA using different stochastic ordering}
For a given copula and an $n$ component parallel system, we look into the monotonicity of the ratio of the distribution functions of parallel systems having dependent  and independent components respectively. We do a similar study for series system. We propose that, we need to analyze the ratio of copula (survival copula) function $C(u)$ ($\hat{C}(u)$) of dependent parallel (series) system and that of independent parallel (series) system $\prod_{i=1}^{n}u_{i}(t)$ ($\prod_{i=1}^{n}\hat{u}_{i}(t)$), where ${u}=({u}_1(t),{u}_2(t),\ldots,{u}_n(t)),$  $t>0.$ \\

We state the following propositions without proof and apply them in obtaining the results in Table~ \ref{tab1}. The following theorem speaks about existence of $rhr$ ordering among two parallel systems having dependent and independent components respectively, under certain conditions. The other weaker orderings follow as usual. As a result, we are able to infer about over-assessment and under assessment in the corresponding error functions.
\begin{t1}
\label{thp}
If $\frac{C(u)}{\Big(\prod_{i=1}^{n}u_{i}\Big)}$ is increasing (decreasing) in $t$ then  $T_{P}^{D} \geq_{rhr} (\leq_{rhr})T_{P}^{I}.$ This in turn gives  $T_{P}^{D} \geq_{st} (\leq_{st}) T_{P}^{I},$ i.e., there is  UA (OA) in reversed hazard rate and in survival function.\end{t1}
{\bf Proof.} From (\ref{eq2p}), one can see that the distribution function of a parallel system having $n$ dependent components is $F^{D}_{P}(t)=C(u(t)),$ whereas for a parallel system having $n$ independent components, the corresponding distribution function is $F^{I}_{P}(t)=\prod_{i=1}^{n}u_{i}(t)$ for $t\geq 0.$ So, the monotonicity of the ratio, viz.,  $$\frac{F^{D}_{P}(t)}{F^{I}_{P}(t)}=\frac{C(u)}{\Big(\prod_{i=1}^{n}u_{i}\Big)}$$ establishes the resultant $rhr$ ordering between the lifetime random variables $T^{D}_{P}$ and $T^{I}_{P}$ respectively. The rest of the lines of proof follow immediately.$\hfill\Box$\\

${}$\hspace{0.5cm}The following theorem speaks about existence of $hr$ ordering among two series systems having dependent and independent components respectively, under certain conditions. A similar exercise as in Theorem \ref{thp}, related to OA (UA) is done for series systems.
\begin{t1}
\label{ths}
Let $\hat{C},$ denote a survival copula with marginal survival functions $\bar{F}_{i}$ for $i=1,2,\ldots,n.$ If $\frac{\hat{C}(u)}{\Big(\prod_{i=1}^{n}\hat{u}_{i}\Big)}$ is increasing (decreasing) in $t$ then $T_{S}^{D} \geq_{hr} (\leq_{hr}) T_{S}^{I}.$ Consequently,   $T_{S}^{D} \geq_{mrl} (\leq_{mrl}) T_{S}^{I}$ and $T_{S}^{D} \geq_{st} (\leq_{st}) T_{S}^{I},$ i.e., there is  UA (OA) in hazard rate, mean residual life and survival function.
\end{t1}
{\bf Proof.} Using (\ref{eq1s}) and similar justification as in Theorem \ref{thp} , we arrive at the result.  $\hfill\Box$ \\
${}$\hspace{0.5cm}We now state a theorem without proof that helps us to establish stochastic ordering between lifetimes of a series (parallel) system with dependent and independent components if the corresponding ordering between parallel (series) system is known, provided dependent components follow a radially symmetric copula. 
\begin{t1}
For a radially symmetric copula,  $\frac{C(u(t))}{\Big(\prod_{i=1}^{n}u_{i}(t)\Big)}$ increases (decreases) with $t$ if and only if $\frac{\hat{C}(\hat{u}(t))}{\Big(\prod_{i=1}^{n}\hat{u}_{i}(t)\Big)}$ is decreases (increases) with $t$.$\hfill\Box$
\end{t1}
{\bf Proof.}
\section{Main Results}
In this section, we deal with several copulas available in literature and find its application in studying error assessment, particularly over and under assessment arising due to the one of the biased assumption that components work independently, though they are not, so as to make the computations easy and simple. \\
${}$\hspace{0.5cm}All through out the upcoming sections, we use  notations viz., $C_1, C_2,\ldots,C_{11}$  so as to comprehend the present study involving several copulas. In each case, we denote the lifetime of a parallel and a series system with components being dependently working by $T^{D}_P$ and $T^{D}_S$ respectively independent of the fact that which particular copula they follow. Similarly, we use the notations $T^{I}_P$ and $T^{I}_S$ respectively if the components of the corresponding system are working independently. In every subsection, we stick to the said notations so as to ease out the style of writing this article. 
\subsection{Independent Copula}
We denote the independent copula by $C_1$ where $C_{1}(u)=\prod_{i}^{n} u_{i}.$ In error assessment, $C_1$ comes as a comparator with every other copula under consideration. 
 \subsection{Farlie-Gumbel-Morgenstern (FGM)}
A commonly known copula is Farlie-Gumbel-Morgenstern (FGM) copula, once proposed by Eyraud (1936)-Farlei (1960)-Gumbel (1958, 1960)-Morgenstern (1956). Various forms of the said copula exist in literature. The FGM copula is helpful in modelling of directional dependence in the field of genetic engineering (Kim et al. (2009)), exchange markets (Jung et al. (2008)) etc. Here, we focus on a  $n$ dimensional FGM copula, denoted by $C_2(u)$ as given
\begin{equation}
\label{fgm1}
C_2(u)=\Big(\prod_{i}^{n} u_{i}\Big) \Big\{1+\alpha\prod_{i=1}^{n}\Big(1-u_{i}\Big)\Big\},\alpha\in [-1,1].
\end{equation}
Since, FGM copula is radially symmetric, the corresponding survival copula can be obtained by replacing $u$ by $1-u$ in the copula function, i.e., we have
\begin{equation}
\label{fgm2}
\hat{C}_2(\hat{u})=\Big(\prod_{i}^{n} \hat{u}_{i}\Big) \Big\{1+\alpha\prod_{i=1}^{n}\Big(1-\hat{u}_{i}\Big)\Big\},\alpha\in [-1,1].\end{equation}
${}$\hspace{0.8cm}In the next theorem, we establish that there exists RHR  ordering between lifetime of $n$ component parallel systems having dependent and independent components respectively if they follow FGM copula (dependent and independent counterpart).  As a result, we are able to throw light on other orderings, $st$ ordering for parallel systems. Using this, we claim that $hr$ ordering exists between lifetime of $n$ component series systems with dependent and independent components respectively following FGM copula.   The proof is easy to see and hence omitted.
\begin{t1}
If $C_2$ and $C_1$ refer to FGM and independent copula respectively then the following hold: \begin{enumerate}
    \item[(i)]If $\alpha$ is positive (negative) then $\frac{C_2}{C_1}$ is decreasing (increasing) in $t.$ \\${}$\hspace{0.4cm}Equivalently, if  $\alpha$ is positive (negative) then $T_{P}^{D}\leq_{rhr} (\geq_{rhr}) ~T_{P}^{I}.$  This gives if $\alpha$ is positive (negative) then $T_{P}^{D}\leq_{st} (\geq_{st})T_{P}^{I}.$
    \item[(ii)]If $\alpha$ is positive (negative) then  $\frac{\hat{C}_2}{\hat{C}_1}$ is increasing (decreasing) in $t.$ \\${}$\hspace{0.4cm}Equivalently if $\alpha$ is positive (negative) then  $T_{S}^{D}\geq_{hr} (\leq_{hr}) T_{S}^{I}$, thereby giving $T_{S}^{D}\geq_{mrl} (\leq_{mrl}) T_{S}^{I}$, and  $T_{S}^{D}\geq_{st} (\leq_{st}) T_{S}^{I}$ according as  $\alpha$ is positive (negative).$\hfill\Box$
\end{enumerate} 
\end{t1}
${}$\hspace{0.8cm}The next theorem proves that there is hazard rate ordering between parallel and series system each having dependent components following FGM copula.
\begin{t1}
\label{fgm}
If $C_2$ and $\hat{C}_2$ denote the copula and survival copula of $FGM$ then $\frac{C_2}{\hat{C}_2}$ is increasing in $t$ for any $\alpha\in[-1,1].$ This gives $T_{P}^{D}\geq_{hr}T_{S}^{D}$ if components follow FGM copula.
\end{t1}
{\bf Proof.} Here, 
\begin{eqnarray}
\frac{C_2}{\hat{C}_2}&=&\frac{\Big(\prod_{i}^{n}{u}_{i}\Big)\Big(1+\alpha \prod_{i}^{n} \hat{u}_{i}\Big)}{\Big(\prod_{i}^{n}\hat{u}_{i}\Big)\Big(1+\alpha \prod_{i}^{n} {u}_{i}\Big)}\nonumber\\
\label{fgme}
&=&\Big\{\alpha+\Big(\frac{1}{\prod_{i=1}^{n}\hat{u}_{i}}\Big)\Big\}\Big\{\frac{1}{\alpha+\Big(\frac{1}{\prod_{i=1}^{n}u_i}\Big)}\Big\}
\end{eqnarray}
which is increasing in $t$ if $\alpha\in [-1,1].$ One shall simultaneously note that each of the factors in right hand side of (\ref{fgme}) are non-negative even if $\alpha\in [-1,0].$ This completes the proof.$\hfill\Box$\\

${}$\hspace{0.8cm}In the following example, we take up Gumbel-II distribution, and note that it can be obtained from FGM survival copula for a particular choice of $\hat{u}_i.$
\begin{ex}
The survival function of bivariate Gumbel II is
\begin{equation}
\label{ex2}
  \bar{F}(t_1,t_2)=\Big\{1+\alpha\left(1-e^{-\lambda_1t_1}\right)\left(1-e^{-\lambda_2t_2}\right)\Big\}e^{-\lambda_1t_1-\lambda_2t_2},\quad |\alpha|<1,
\lambda_i,t_i>0,\;i=1,2. 
\end{equation}
The survival function of Gumbel-II can be obtained from the survival copula of Farlie Gumbel  Morgenstern by taking $\hat{u}_{i}(t_i)=\exp(-\lambda_{i}t_i),t \geq 0,$ and $i=1,2.$
To validate our point, we see that (\ref{ex2}) can be written as  $$\bar{F}(x_1,x_2)=\Big\{\prod_{i=1}^{2}\hat{u}_{i}(t_i)\Big\}\Big\{1+\alpha \prod_{i=1}^{2}(1-\hat{u}_{i}(t_i))\Big\},$$ with $\hat{u}_{i}(t)=\exp(-\lambda_{i}t_i).$ $\hfill\Box$
\end{ex}
\subsection{Fischer and K{\"o}ck's (2012) Copula}
 There are different forms of Fischer and K{\"o}ck's (2012) copula. Here, we represent it by $C_3(u)$ where $$C_3(u)=\Big(\prod_{i}^{n} u_{i}\Big) \Big\{1+\alpha \prod_{i=1}^{n}\Big(1-u_{i}^{1/r}\Big)\Big\}^{r}, r\geq 1; -1\leq \alpha \leq 1.$$ 
Notably, Fischer and K{\"o}ck's copula is a generalization of FGM copula for $r=1.$ Since, the said copula is radially symmetric, we write the corresponding survival copula as $$\hat{C}_3(u)=\Big(\prod_{i}^{n} \hat{u}_{i}\Big) \Big\{1+\alpha \prod_{i=1}^{n}\Big(1-\hat{u}_{i}^{1/r}\Big)\Big\}^{r},r\geq 1; -1\leq \alpha \leq 1.$$ 
  \begin{t1}
\label{fgm}
If $C_3$ and $\hat{C}_3$ denote the copula and survival copula proposed by Fischer and K{\"o}ck's then the following hold:
\begin{enumerate}
   \item[(i)]If $\alpha$ is positive (negative) then $\frac{C_3}{C_1}$ is decreasing (increasing) in $t.$ \\${}$\hspace{0.4cm}Equivalently, if  $\alpha$ is positive (negative) then $T_{P}^{D}\leq_{rhr} (\geq_{rhr}) ~T_{P}^{I}.$  This gives if $\alpha$ is positive (negative) then $T_{P}^{D}\leq_{st} (\geq_{st})T_{P}^{I}.$
    \item[(ii)]If $\alpha$ is positive (negative) then  $\frac{\hat{C}_3}{\hat{C}_1}$ is increasing (decreasing) in $t.$\\${}$\hspace{0.4cm} Equivalently if $\alpha$ is positive (negative) then  $T_{S}^{D}\geq_{hr} (\leq_{hr}) T_{S}^{I}$, thereby giving $T_{S}^{D}\geq_{mrl} (\leq_{mrl}) T_{S}^{I}$, and  $T_{S}^{D}\geq_{st} (\leq_{st}) T_{S}^{I}$ according as  $\alpha$ is positive (negative).$\hfill\Box$
\end{enumerate}
\end{t1}
 \begin{t1}
\label{fgm}
If $C_3$ and $\hat{C}_3$ denote the copula and survival copula of Fischer and K{\"o}ck's copula then $\frac{C_3}{\hat{C}_3}$ is increasing in $t$ for any $\alpha\in[-1,1].$ This gives $T_{P}^{D}\geq_{hr}T_{S}^{D}$ if components follow Fischer and K{\"o}ck's copula.
\end{t1}
{\bf Proof.} We note that,
\begin{eqnarray}
\frac{C_3}{\hat{C}_3}&=&\frac{\Big(\prod_{i}^{n} {u}_{i}\Big) \Big\{1+\alpha \prod_{i=1}^{n}\Big(1-{u}_{i}^{1/r}\Big)\Big\}^{r}}{\Big(\prod_{i}^{n} \hat{u}_{i}\Big) \Big\{1+\alpha \prod_{i=1}^{n}\Big(1-\hat{u}_{i}^{1/r}\Big)\Big\}^{r}}\nonumber\\
\label{fgme}
\label{ff}
&=&\Big[\frac{1+\alpha \prod_{i=1}^{n} (1-u_i^{1/r})}{\prod_{i=1}^{n}(1-u_i)^{\frac{1}{r}}}\Big]^{r}\Big[\frac{1+\alpha \prod_{i=1}^{n} (1-\hat{u}_i^{1/r})}{\prod_{i=1}^{n}(1-\hat{u}_i)^{\frac{1}{r}}}\Big]^{-r}\nonumber\\
&=&XY (say), 
\end{eqnarray}
where $X$ and $Y$ refer to the two consecutive products as given in (\ref{ff}).
We first study the monotonic behaviour of $V=\frac{1-u_i^{\frac{1}{r}}}{(1-u_i)^{\frac{1}{r}}}$ and $W=\frac{1-\hat{u}_i^{\frac{1}{r}}}{(1-\hat{u}_i)^{\frac{1}{r}}}$ respectively  to conclude about $X$ and $Y$ respectively. As a simple exercise, one may check that 
\begin{eqnarray}
 \ln V&=&\ln ((1-u_i)^{\frac{1}{r}})-\frac{1}{r}\ln(1-u_i)\nonumber
 \end{eqnarray}
 After differentiating, we get 
 \begin{eqnarray}
\frac{1}{V}\frac{dV}{dt}=\frac{1}{r}\Big\{\frac{1-u_i^{{\frac{1}{r}-1}}}{(1-u_i)(1-u_i^{\frac{1}{r}}
  )}\Big\}\frac{du_i}{dt}\nonumber  
 \end{eqnarray} 
which is negative   
 draw the analogous conclusion for $N=\frac{1+\alpha \prod_{i=1}^{n} (1-\hat{u}_i^{1/r})}{\prod_{i=1}^{n}(1-\hat{u}_i)^{\frac{1}{r}}}.$ which is increasing in $t$ if $\alpha\in [-1,1].$ One shall simultaneously note that each of the factors in right hand side of (\ref{fgme}) are non-negative even if $\alpha\in [-1,0].$ This completes the proof.$\hfill\Box$\\
\subsection{Gumbel-Hougaard}
Here, we denote the copula, called as Gumbel-Hougaard by $C_5(u)=\exp\Big[\Big\{-\sum_{i=1}^{n}\big(-\ln u_i\big)^{\alpha}\Big\}^{\frac{1}{\alpha}}\Big].$
\begin{t1}
For Gumbel-Hougaard copula, denoted by $C_5$ the following conditions hold:
\begin{equation}
\frac{C_5(v)}{C_1(v)}=\frac{\exp\Big[-\Big\{\sum_{i=1}^{n}\Big(-\ln v_i\Big)^{\alpha}\Big\}^{\frac{1}{\alpha}}\Big]}{\exp\Big[\ln (\prod_{i=1}^{n}v_{i})\Big]}, \alpha \geq 1\nonumber
\end{equation}
is decreasing  in $t\geq 0,$ if $v=(v_1,v_2,\ldots,v_n)$ where $v_i=u_i=F_{i}(t)$  and is increasing in $t$, if $v=(v_1,v_2,\ldots,v_n)$ where $v_i=\hat{u}_i=\hat{F}_{i}(t).$ 
\end{t1}
{\bf Proof.} We note that
\begin{equation}
\frac{C_5(v)}{C_1(v)}=\frac{C_5(v)}{\prod_{i=1}^{n}v_{i}}=\frac{\exp\Big[-\Big\{\sum_{i=1}^{n}\Big(-\ln v_i\Big)^{\alpha}\Big\}^{\frac{1}{\alpha}}\Big]}{\exp\Big[\ln (\prod_{i=1}^{n}v_{i})\Big]},~ \alpha \geq 1.\nonumber
\end{equation}
Taking $v=u$, we find that $ \frac{C_5(u)}{\prod_{i=1}^{n}u_{i}}=\exp\big(-\phi(t)\big),$ where \begin{equation}
\phi(t)=\Big\{\sum_{i=1}^{n}\Big(-\ln u_i\Big)^{\theta}\Big\}^{\frac{1}{\theta}}-\sum_{i=1}^{n}\ln u_i =\Big(\sum_{i=1}^{n}x_i^{\alpha}\Big)^{\frac{1}{\alpha}}-\sum_{i=1}^{n}x_i,\nonumber
\end{equation}
with $x_i=-\ln u_{i}\geq 0.$
We note that
\begin{eqnarray}
\frac{d}{dt}\phi(t)&=&\sum_{i=1}^{n}x_{i}^{\alpha-1}\Big(\sum_{i=1}^{n}x_{i}^{\alpha}\Big)^{\frac{1}{\alpha}-1}\frac{dx_i}{dt}\nonumber\\
&=&\sum_{i=1}^{n}\Bigg[\Bigg\{\Bigg(1+\Big(\frac{x_1}{x_i}\Big)^{\alpha}+\ldots+\Big(\frac{x_{i-1}}{x_i}\Big)^{\alpha}+\Big(\frac{x_{i+1}}{x_{i}}\Big)^{\alpha}+\ldots+\Big(\frac{x_{n}}{x_i}\Big)^{\alpha}\Bigg)^{\frac{1}{\alpha}-1}-1\Bigg\}\frac{dx_i}{dt}\Bigg]\geq 0,\nonumber
\end{eqnarray}
since $\alpha_i\geq 1$ and $\frac{dx_i}{dt}\leq 0$ for all $i.$ Similarly considering $v=\hat{u}$ we arrive at the conclusion. This proves the theorem.$\hfill\Box$\\
The following remark helps us to establish reversed hazard rate (hazard rate) ordering between lifetimes of parallel (series) system having dependent and independent components. It also throws light on other orderings that follow from a stronger ordering.
\begin{r1}
If $C_5$ represents the distribution copula of Gumbel-Hougaard copula then $T_{P}^{D}\leq_{rhr} T_{P}^{I}$ which in turn implies  $T_{P}^{D}\leq_{st} T_{P}^{I}.$ However, if $C_5$ represents the survival copula of Gumbel-Hougaard copula then $T_{S}^{D}\geq_{hr} T_{S}^{I}$ giving $T_{S}^{D}\geq_{mrl} T_{S}^{I}$ and $T_{S}^{D}\geq_{st} T_{S}^{I}.$$\hfill\Box$
\end{r1}
Below, we give an example of a bivariate distribution that can be obtained from Gumbel-Hougaard survival copula.
\begin{ex}
The survival function of Gumbel-III given by $$\bar{F}(t_1,t_2)=e^{-((\lambda_{1} t_{1})^{\alpha} + ( \lambda_{2} t_{2})^{\alpha} )^{\frac{1}{\alpha}}}, \quad
\lambda_i,t_i>0,\;i=1,2;~\alpha \geq 1; $$ can be obtained from  survival copula of Gumbel-Hougaard by taking $\hat{u}_{i}(t_i)=\exp(-\lambda_{i}t_i),t_i \geq 0,$ for $i=1,2$
and write as
\begin{eqnarray}
\bar{F}(t_1,t_2)=\exp\Big(-\sum_{i=1}^{2}(-\ln \hat{u}_{i}(t_i))^{\alpha}\Big)^{\frac{1}{\alpha}}\nonumber
\end{eqnarray}
since $\lambda_it_i=-\ln \hat{u}_i(t_i).$   $\hfill\Box$
\end{ex}
\subsection{  Gumbel-Barnet   }
In this subsection, we take up Gumbel-Barnet copula. We denote the corresponding copula by $C_6(\cdot)$ given by
$C_6(u)=\Big(\prod_{i=1}^{n}u_{i}\Big)\exp\Big(-\alpha \prod_{i=1}^{n}\ln u_i\Big).$\\
We state the following theorem without proof.
\begin{t1}
For Gumbel-Barnet, denoted by $C_{6}, \frac{C_6(v)}{C_1(v)}$ is decreasing (increasing) in $t$ if  $v=u$ ($v=\hat{u}$).
\end{t1}
As a result, we write the following remark to give its importance in establishing stochastic orderings between lifetime of parallel and series systems.
\begin{r1}
If $C_5$ represents the distribution copula of Gumbel-Barnet copula then i.e., $T_{P}^{D}\leq_{rhr} T_{P}^{I},$ giving $T_{P}^{D}\leq_{st} T_{P}^{I}.$ On the other hand, if $C_5$ represents the survival copula of Gumbel-Barnet copula then $T_{S}^{D}\geq_{hr} T_{S}^{I}.$ The ordering weaker than $hr$ order follows immediately.
\end{r1}
We now give an example to illustrate the use of Gumbel-Barnet survival copula in obtaining the survival function of Gumbel-I bivariate exponential distribution.
\begin{ex}
The joint survival function of bivariate Gumbel I exponential distribution is
$$\bar{F}(x_1,x_2)=e^{-\lambda_1x_1-\lambda_2x_2-\lambda_{12}x_1x_2},\quad
\lambda_{12}\ge 0; ~0 \leq \lambda_{12} \leq \lambda_{1} \lambda_{2}, $$
$\lambda_i,x_i>0,\;i=1,2;$\\
Taking $\hat{u}_{i}(x_i)=\exp(-\lambda_{i}x_i),$ for $i=1,2$ we find that
\begin{eqnarray}
\bar{F}(x_1,x_2)&=&e^{-\lambda_1x_1}e^{-\lambda_2x_2}e^{-\lambda_{12}x_1x_2}\nonumber\\
&=&\hat{u_1}(x_1)\hat{u_2}(x_2)\exp\Big\{\Big(-\frac{\lambda_{12}}{\lambda_1\lambda_2}\Big)\big(-\lambda_1 x_1\big)\big(-\lambda_2 x_2\big)\Big\}\nonumber\\
&=&\hat{u_1}(x_1)\hat{u_2}(x_2)\exp\Big\{\Big(-\frac{\lambda_{12}}{\lambda_1\lambda_2}\Big)\big(\ln \hat{u_1}(x_1) \ln \hat{u_2}(x_2) \big)\Big\}\nonumber\\
&=&\big(\prod_{i=1}^{2}\hat{u}_{i}\big)\exp\big(-\alpha \prod_{i=1}^{2}\ln \hat{u}_i\big)\nonumber
\end{eqnarray}
where $\alpha=\frac{\lambda_{12}}{\lambda_1\lambda_2}$ and this is a particular  realization of Gumbel-Barnet copula. 
\end{ex}
\subsection{Nelson Ten}	
 $C_7(u)=\frac{\prod_{i=1}^{n}u_{i}}{\Big\{1+\prod_{i=1}^{n}(1-u_i^{\alpha})\Big\}^{\frac{1}{\alpha}}}$
\begin{t1}
 $\frac{C_7}{C_1}\uparrow,$ i.e.,$T_{P}^{D}\geq_{rhr} T_{P}^{I}\Rightarrow T_{P}^{D}\geq_{st} T_{P}^{I}$.
			 $\frac{\hat{C}_7}{\hat{C}_1}\downarrow,$ i.e.,$T_{S}^{D}\leq_{hr} T_{S}^{I}$ giving $T_{S}^{D}\leq_{mrl} T_{S}^{I}$ and $T_{S}^{D}\leq_{st} T_{S}^{I}$
\end{t1}

\subsection{    Marshall and Olkin}
We know that Marshall-Olkin copula, here denoted by $C_8$ is given by $$C_8(u)=\big(\prod_{i=1}^{n}u_{i}\big)\min\big\{u_{1}^{\alpha_1}, u_{2}^{\alpha_2},\ldots, u_{n}^{\alpha_n}\big\},$$ where $\alpha_i >0$ for $i=1,2,\ldots,n.$
Marshall-Olkin copula being radially symmetric, Marshall-Olkin survival copula here denoted by $\hat{C}_8$ is given by $$\hat{C}_8(u)=\big(\prod_{i=1}^{n}\hat{u}_{i}\big)\min\big\{\hat{u}_{1}^{\alpha_1}, \hat{u}_{2}^{\alpha_2},\ldots, \hat{u}_{n}^{\alpha_n}\big\}.$$
Before proceeding for further study required for error analysis, in the following example, we mention the joint survival function of  Marshall and Olkin exponential distribution. We show that the joint survival function can be expressed in the form of Marshall and Olkin survival copula with a particular choice of marginal survival functions, i.e., $\hat{u}_i's.$
\begin{ex}
 Marshall and Olkin Bivariate Exponential Distribution:
$$\bar{F}(x_1,x_2)=e^{-\lambda_1x_1-\lambda_2x_2-\lambda_{12}\max(x_1,x_2)},\quad
\lambda_i,x_i>0,\;i=1,2;\;\lambda_{12}> 0;$$
\begin{eqnarray}
    \bar{F}(x_1,x_2)&=&e^{-\lambda_1x_1}e^{-\lambda_2x_2}e^{-\lambda_{12}\max{(x_1,x_2)} }\nonumber\\
&=& e^{-\lambda_1x_1}e^{-\lambda_2x_2}e^{-\max{(\lambda_{12}x_1,\lambda_{12}x_2)} }\nonumber\\
&=&e^{-\lambda_1x_1}e^{-\lambda_2x_2}e^{\min{(-\lambda_{12}x_1,-\lambda_{12}x_2}) }\nonumber\\
&=& e^{-\lambda_1x_1}e^{-\lambda_2x_2}\min\{e^{{-\lambda_{12}x_1} } ,e^{-\lambda_{12}x_2}\}\nonumber\\
&=& e^{-\lambda_1x_1}e^{-\lambda_2x_2}\min\{e^{{\frac{\lambda_{12}}{\lambda_{1}}(-\lambda_{1}x_1)} } ,e^{{\frac{\lambda_{12}}{\lambda_{2}}(-\lambda_{2}x_2)} } \}\nonumber\\
&=&\hat{u}_{1}(x_1)\hat{u}_{2}(x_2)\min\{\hat{u}^{\alpha}_{1}(x_1),\hat{u}^{\beta}_{2}(x_2)\}\nonumber
\end{eqnarray}
where $\hat{u}_{i}=e^{-\lambda_i x_{i}},$ for $i=1,2,$  $\alpha=\frac{\lambda_{12}}{\lambda_1},$ and $\beta=\frac{\lambda_{12}}{\lambda_2}.$$\hfill\Box$
\end{ex}
In the following example, we mention the joint survival function of Block and Basu exponential distributions. We show that the said distribution, joint survival function can be expressed in the form of Marshall and Olkin survival copula with a particular choice of marginal survival functions, i.e., $\hat{u}_i's.$
\begin{ex}
  Block and Basu Bivariate Exponential Distribution:
$$\bar{F}_7(x_1,x_2)=\frac{\lambda^*}{\lambda}e^{-\lambda_1x_1-\lambda_2x_2-\lambda_{12}\max(x_1,x_2)}-\frac{\lambda_{12}}{\lambda}e^{-\lambda^*\max(x_1,x_2)},\quad $$
where $\lambda_i,x_i>0,\;i=1,2,\lambda_{12}\ge 0,\lambda=\lambda_1+\lambda_2,\lambda^*=\lambda+\lambda_{12};$ Let $\bar{U}_{i}=e^{-\lambda_i}x_i$ for $i=1,2.$
\\

Let $\theta=\frac{\lambda^*}{\lambda}$ then $1-\theta=-\frac{\lambda_{12}}{\lambda}.$
\begin{eqnarray}
    \bar{F}(x_1,x_2)&=&\theta e^{-\lambda_1x_1}e^{-\lambda_2x_2}e^{-\lambda_{12}\max{(x_1,x_2)} }+(1-\theta)e^{-\lambda^*\max{(x_1,x_2)} }\nonumber\\
&=& \theta e^{-\lambda_1x_1}e^{-\lambda_2x_2}e^{-\max{(\lambda_{12}x_1,\lambda_{12}x_2)} }+(1-\theta)e^{-\max{(\lambda^*x_1,\lambda^*x_2)} }\nonumber\\
&=&\theta e^{-\lambda_1x_1}e^{-\lambda_2x_2}e^{\min{(-\lambda_{12}x_1,-\lambda_{12}x_2}) }+(1-\theta)e^{\min{(-\lambda^*x_1,-\lambda^*x_2)} }\nonumber\\
&=& \theta e^{-\lambda_1x_1}e^{-\lambda_2x_2}\min\{e^{{-\lambda_{12}x_1} } ,e^{-\lambda_{12}x_2}\}+(1-\theta)\min\{e^{{-\lambda^*x_1} } ,e^{-\lambda^*x_2}\}\nonumber\\
&=&  \theta e^{-\lambda_1x_1}e^{-\lambda_2x_2}\min\{e^{{\frac{\lambda_{12}}{\lambda_{1}}(-\lambda_{1}x_1)} } ,e^{{\frac{\lambda_{12}}{\lambda_{2}}(-\lambda_{2}x_2)} } \}+(1-\theta)\min\{e^{{\frac{\lambda^*}{\lambda_{1}}(-\lambda_{1}x_1)} } ,e^{{\frac{\lambda^*}{\lambda_{2}}(-\lambda_{2}x_2)} } \}\nonumber\\
&=&  \theta\hat{u}_{1}(x_1)\hat{u}_{2}(x_2)\min\{\hat{u}^{\alpha}_{1}(x_1),\hat{u}^{\beta}_{2}(x_2)\}+(1-\theta)\min\{\hat{u}^{\delta}_{1}(x_1),\hat{u}^{\gamma}_{2}(x_2)\}\nonumber
\end{eqnarray}
where $\alpha=\frac{\lambda_{12}}{\lambda_1},$ $\beta=\frac{\lambda_{12}}{\lambda_2}$, $\delta=\frac{\lambda^*}{\lambda_1},$ $\gamma=\frac{\lambda^*}{\lambda_2}$  and $\hat{u}_{i}=e^{-\lambda_i x_{i}},$ for $i=1,2.$
Here, one can note that $\delta$ and $\gamma$ are related to $\alpha$ and $\beta$ respectively.
$\delta=1+\alpha+\frac{\lambda_2}{\lambda_1}$ and $\gamma=1+\beta+\frac{\lambda_1}{\lambda_2}.$$\hfill\Box$
  \end{ex}
Before we proceed to study the monotonic behavior of the ratio of Marshall and Olkin copula and independent copula, we state and prove two lemmas used to describe some properties of $u_i$ and $\hat{u}_i$.
\begin{l1}
\label{minl1}
 $\min\{u_1(t), u_2(t),\ldots,u_n(t) \}$ is increasing in $t\geq 0.$
\end{l1}
{\bf Proof.}
We first prove that the $\min\{u_1(t), u_2(t) \}$ is increasing in $t\geq 0.$
Let $b\in R$ such that for some $\epsilon>0,$
for $x\in (b-\epsilon, b),$ $\min\{F_1(x),F_2(x)\}=F_2(x)$ and for $x\in (b, b+\epsilon),$ $\min\{F_1(x),F_2(x)\}=F_1(x).$ Since $F_1(z)$ is increasing in $z,$ so for all $x\in [b,b+\epsilon)$ and for all $y\in [b-\epsilon,b),$ $F_1(y)\leq F_1(x).$ This implies $F_2(y)\leq F_1(x)$ since $F_2(y)=\min (F_1(y),F_2(y)),$ i.e., $\min \{F_1(y),F_2(y)\}\leq \min \{F_1(x),F_2(x)\}.$ Additionally, for all $x\leq y \in [b,b+\epsilon),$ $\min\{F_1(x), F_2(x)\}\leq \min\{F_1(y), F_2(y)\}$ since $\min\{F_1(x), F_2(x) \}=F_{1}(x)$ and $\min\{F_1(y), F_2(y) \}=F_{1}(y)$ and $F_1(x)$ is increasing. Similarly, $\min \{F_{1}(x), F_{2}(x)\}$ is increasing in $(b-\epsilon,b).$ Thus, $\min\{F_{1}(x),F_{2}(x)\}$ is increasing in $(b-\epsilon,b+\epsilon).$ We have already established that the minimum of two increasing functions is increasing. Now, consider the increasing functions $f_1, f_2, \ldots, f_n$. Let us define $g_k(x)=\min\{F_1(x),F_2(x),\ldots,F_k(x)\},$ for all $k\in\{2,\ldots,n\}.$ Clearly, $g_2$ is increasing. That implies $g_3(x)=\min\{g_2(x),F_3(x)\}$ is also increasing. In fact, if $g_k$ is increasing, then $g_{k+1}(x)=\min\{g_k(x),F_{k+1}(x)\}$ is also increasing. Hence, by the principle of mathematical induction, $g_n$ is increasing. $\hfill\Box$
\begin{l1}
\label{minl2}
 $\min\{\hat{u}_1(t), \hat{u}_2(t),\ldots,\hat{u}_n(t) \}$ is decreasing in $t\geq 0.$
\end{l1}
Let $u_1$, $u_2$ be decreasing functions, i.e., $\forall x<y$, $$u_1(x)\geq u_1(y) \text{, and }u_2(x)\geq u_2(y).$$
If $a\in\mathbb{R}$ be a point such that for some $\epsilon>0$, $\forall x\in (a-\epsilon,a+\epsilon)$, $\min\{u_1,u_2\}$ is identically equal to either of $u_1$ or $u_2$, then clearly, $\min\{u_1,u_2\}$ is also decreasing in $(a-\epsilon,a+\epsilon)$.
Let $b\in\mathbb{R}$ be a point such that for some $\epsilon>0$,  for all $x\in(b-\epsilon,b)$, $u_1(x)\leq u_2(x),$ and $\forall x\in[b,b+\epsilon)$, $u_1(x)\geq u_2(x).$
Therefore, if $b-\epsilon<x<y<b$, then $\min\{u_1,u_2\}(x)= u_1(x)\geq u_1(y)=\min\{u_1,u_2\}(y).$
	If $b\leq x<y< b+\epsilon$, then $\min\{u_1,u_2\}(x)=u_2(x)\geq u_2(y)=\min\{u_1,u_2\}(y).$
	If $b-\epsilon<x<b\leq y< b+\epsilon$, then $\min\{u_1,u_2\}(x)= u_1(x)\geq u_1(y)\geq u_2(y)=\min\{u_1,u_2\}(y).$
	Therefore, $\min\{u_1,u_2\}$ is a decreasing function.\\
${}$\hspace{0.8cm} We have already established that the minimum of two decreasing functions is decreasing. Now, consider the decreasing functions $u_1, u_2, \ldots, u_n$. Let us define $g_k(x)=\min\{u_1(x),u_2(x),\ldots,u_k(x)\},$for all $k\in\{2,\ldots,n\}.$
	Clearly, $g_2$ is decreasing. That implies $g_3(x)=\min\{g_2(x),f_3(x)\}$ is also decreasing. In fact, if $g_k$ is decreasing, then $g_{k+1}(x)=\min\{g_k(x),f_{k+1}(x)\}$ is also decreasing. Hence, by the principle of mathematical induction, $g_n$ is decreasing.$\hfill\Box$\\
 We now give the following theorem using Lemma \ref{minl1}.
\begin{t1}
\label{rhrm}
 $\frac{C_{8}(u)}{\prod_{i=1}^{n}u_{i}}=\min\{u_i^{\alpha_i}(t):1\leq i\leq n\}$ is increasing in $t\geq 0,$ where $u_{i}(t)=F_i(t).$
  \end{t1}
Using Theorem \ref{rhrm}, in the following remark we establish that there exists reversed hazard rate order and stochastic order among the lifetimes of parallel system comprising of finite number of dependent (following Marshall Olkin copula) and independent components respectively.
 \begin{r1}
If $T_{P}^{D}$ and Since $\frac{{C}_8}{{C}_1}$ is increasing in $t$,it follows that $T_{P}^{D}\geq_{rhr} T_{P}^{I}.$ This implies that $T_{P}^{D}\geq_{st} T_{P}^{I}.$
  \end{r1}
  In the next theorem, we show that
  \begin{t1}
 $\frac{\hat{C}_{8}(u)}{\prod_{i=1}^{n}\hat{u}_{i}}=\min\{\hat{u}_i^{\alpha_i}(t):1\leq i\leq n\}$ is decreasing in $t\geq 0,$ where $\hat{u}_{i}(t)=\bar{F}_i(t).$
 \end{t1}
  \begin{c1}
Since, $\frac{\hat{C}_8}{\hat{C}_1}$ if $\alpha \geq 0$, is decreasing in $t\geq 0,$ it follows that $T_{P}^{D}\leq_{hr} T_{P}^{I} \Rightarrow T_{P}^{D}\geq_{st} T_{P}^{I}$
\end{c1}

{\bf Proof.} $\hfill\Box$
     \subsection{(AMH) Ali-Mikhail-Haq}
         $C_9(u)=(1-\alpha)\Big\{\prod_{i=1}^{n}\Big(\frac{1-\alpha}{u_i}+\alpha\Big)-\alpha\Big\}^{-1}$
\begin{t1}
For  AMH copula defined by
$$\frac{C_{9}(u)}{\prod_{i=1}^{n}u_{i}}=(1-\alpha)\Big\{\prod_{i=1}^{n}\Big(1-\alpha+\alpha u_{i}\Big)-\alpha\prod_{i=1}^{n}u_{i}\big)\Big\}^{-1}$$
is decreasing in $t$ if $\alpha \in [0,1]$ and increasing in $t$ if $\alpha \in [-1,0].$ Further, $\frac{C_9}{C_1}\uparrow (\downarrow)$ if $\alpha \in [-1,0] (\alpha \in [0,1]),$ i.e., $T_{P}^{D}\geq_{rhr} (\leq_{rhr}) T_{P}^{I}\Rightarrow T_{P}^{D}\geq_{st} (\leq_{st})T_{P}^{I}$\\
		     $\frac{\hat{C}_9}{\hat{C}_1}\downarrow (\uparrow)$ if $\alpha \in [-1,0] (\alpha \in [0,1]),$ i.e., $T_{S}^{D}\leq_{hr} (\geq_{hr}) T_{S}^{I},$
    So, $T_{S}^{D}\leq_{mrl} (\geq_{mrl}) T_{S}^{I}$,  $T_{S}^{D}\leq_{st} (\geq_{st}) T_{S}^{I}$
\end{t1}
{\bf Proof.} Here,
$$C_{9}(u)=(1-\alpha)\Big\{\prod_{i=1}^{n}\Big(\frac{1-\alpha}{u_i}+\alpha\Big)-\alpha\Big\}^{-1}=(1-\alpha)X^{-1}, -1\leq \alpha \leq 1.$$
where,
$$X^{-1}=\Big\{\prod_{i=1}^{n}\Big(\frac{1-\alpha}{u_i}+\alpha\Big)-\alpha\Big\}^{-1}= \Big(\prod_{i=1}^{n}u_{i}\Big)\Big\{\prod_{i=1}^{n}\Big(1-\alpha+\alpha u_{i}\Big)-\alpha\prod_{i=1}^{n}u_{i}\big)\Big\}^{-1},$$
This gives $\frac{C(u)}{\prod_{i=1}^{n}u_{i}}=(1-\alpha)(A+B)^{-1},$ where $A=\prod_{i=1}^{n}\Big(1-\alpha+\alpha u_{i}\Big),$ and $B=-\alpha\prod_{i=1}^{n}u_{i}.$
We note that if $\alpha \in [0,1],$ then $A\geq 0,$ $B\leq 0,$ $A$ is increasing in $t$ and $B$ is decreasing in $t.$ On the other hand, if $\alpha \in [-1,0],$ then $A,B \geq 0,$ $A$ is decreasing in $t$ and $B$ is increasing in $t.$ So, we need to investigate further to know about the monotonic behavaiour of $A+B$ as function of $t$ which in turn will describe the same about $\frac{C(u)}{\prod_{i=1}^{n}u_{i}}.$ Clearly, \begin{eqnarray}
\frac{d}{dt} (A+B)&=&A \sum_{i=1}^{n}\Big\{\Big(\frac{\alpha}{1-\alpha+\alpha u_{i}}\Big)\Big(\frac{du_{i}}{dt}\Big)\Big\}+ B\sum_{i=1}^{n}\Big(\frac{1}{u_{i}}\frac{du_{i}}{dt}\Big)\nonumber\\
&=&\Big( \prod_{i=1}^{n}\big(1-\alpha+\alpha u_{i}\big)\Big)\sum_{i=1}^{n}\Big\{\Big(\frac{\alpha}{1-\alpha+\alpha u_{i}}\Big)\Big(\frac{du_{i}}{dt}\Big)\Big\}-\alpha\Big(\prod_{i=1}^{n}u_{i}\Big)\sum_{i=1}^{n}\Big(\frac{1}{u_{i}}\frac{du_{i}}{dt}\Big)\nonumber\\
&=&\sum_{i=1}^{n}\Big\{\alpha \Big(\mathop{\prod_{j=1}^{n}}_{j\neq i}\big(1-\alpha+\alpha u_{j}\big)-\mathop{\prod_{j=1}^{n}}_{j\neq i}u_{j}\Big)\frac{du_{i}}{dt}\Big\}, \nonumber
\end{eqnarray}
which is positive if $\alpha \in [0,1]$ and negative if $\alpha \in [-1,0]$, since $1-\alpha+\alpha u_{j} \geq u_{j}$ and $1-\alpha+\alpha u_{j} \geq 0,$ for $\alpha \in [-1,1],$ $j=1,2,\ldots,n.$ Thus, $\frac{C(u)}{\prod_{i=1}^{n}u_{i}}$ is increasing in $t$ if $\alpha \in [-1,0],$ and decreasing in $t$ if $\alpha \in [0,1].$$\hfill\Box$
\subsection{Fischer and Hinzmann's copulas}
    $C_{10}(u)=\Big[\Big\{\alpha\min_{i}u_i\Big\}^{m}+\Big\{(1-\alpha) \prod_{i=1}^{n}u_i\Big\}^{m}\Big]^{1/m}$
    \begin{t1}
for $0\leq \theta \leq 1,$ $\frac{\hat{C}_{10}}{\hat{C}_1}\downarrow,$ i.e.,
 $T_{P}^{D}\leq_{rhr} T_{P}^{I};$ \\
    for $-1\leq \theta \leq 0,$ $\frac{\hat{C}_{10}}{\hat{C}_1}$ is constant ($T_{P}^{D}$ and  $T_{P}^{I}$ have PRHR) until $u_1+u_2<1,$$\hfill\Box$
\end{t1}
     \begin{t1}
$\frac{C_{10}(u)}{C_{1}(u)}$ is decreasing in $t,$ since $\frac{\min_{i}u_{i}(t)}{\prod_{i=1}^{n}u_{i}(t)}$ is decreasing in $t.$ Further, $\frac{C_{10}}{C_1}\downarrow,$ i.e.,$T_{P}^{D}\leq_{rhr} T_{P}^{I}\Rightarrow T_{P}^{D}\leq_{st} T_{P}^{I}$\\
				 $\frac{\hat{C}_{10}}{\hat{C}_1}\uparrow,$ i.e.,$T_{S}^{D}\geq_{hr} T_{S}^{I}$ giving $T_{S}^{D}\geq_{mrl} T_{S}^{I}$ and $T_{S}^{D}\geq_{st} T_{S}^{I}$\\
\end{t1}
{\bf Proof.}
Let $\min_{i}u_{i}(t)=u_{k}(t)$ for all $t \in (a-h,a]=A,$ (say) and $\min_{i}u_{i}(t)=u_{l}(t)$ for all $t \in [a,a+h)=B,$ (say). Then for all $t\in A,$ $V(t)=\frac{\min_{i} u_{i}(t)}{\prod_{i=1}^{n}u_{i}(t)}=\frac{1}{\prod_{ i\neq k}u_{i}(t)}$ and for all $t\in B,$ $V(t)=\frac{1}{\prod_{ i\neq l}^{n}u_{i}(t)}.$ Let $a-h<x<y<a,$ $u_{i}(x)\leq u_{i}(y)$ for all $i \in \{1,2,\ldots,n\},$ which implies $\frac{1}{\prod_{i\neq k}u_{i}(x)}\geq \frac{1}{\prod_{i=1, i\neq k}^{n}u_{i}(x)},$ i.e., $V(x)\geq V(y).$ Similarly, if $a\leq x <y <a+h,$ $\frac{1}{\prod_{i\neq l}u_{i}(x)}\geq \frac{1}{\prod_{i\neq l}u_{i}(x)},$ i.e., $V(x)\geq V(y).$
Now, let $x\in A, y\in B,$ $u_{i}(x)\leq u_{i}(y)$ for all $i \in\{1,2,\ldots,n\}$
giving
$\prod_{ i\neq \{k,l\}}u_{i}(x) \leq \prod_{i\neq \{k,l\}}u_{i}(y).$
Also, 
$u_{k}(x)\leq u_{l}(x)\leq u_{l}(x)\leq u_{k}(y).$
We get $\prod_{ i\neq \{k\}}u_{i}(x) \leq \prod_{i\neq \{l\}}u_{i}(x),$ i.e., $V(x)\geq V(y).$ This completes the proof.$\hfill\Box$\\
        $C_{11}(u)=\prod_{i=1}^{n}u_{i}(t)+\alpha \Big(\prod_{i=1}^{n}u_{i}^{a_i}\Big)\Big( \prod_{i=1}^{n}(1-u_{i})^{b_i}\Big)$
\subsection{Rodriguez -Lellena and Ubeda Flores (2004): An extension}
Rodriguez -Lellena and Ubeda Flores (2004) have introduced a bivariate copula $$C(u_1,u_2)=u_1u_2+\alpha u_1^{a}u_2^{b} (1-u_1)^{c}(1-u_2)^{d},$$ where $a,b,c, d\geq 1.$ We give a simple generalization of the said bivariate copula by defining the corresponding $n$ variate copula $$C_{11}(u)=\prod_{i=1}^{n}u_{i}+\alpha\prod_{i=1}^{n}\Big(u_{i}^{a_i}(1-u_{i})^{b_i}\Big),$$ where $u=(u_1,u_2,\ldots,u_n),$ $a_i,b_i\geq 1,$ and $0\leq\alpha \leq 1.$ We term $C_{11}$ as extension of Rodriguez-Lallena and Ubeda Flores.\begin{t1} For  extended Rodriguez-Lallena and Ubeda Flores copula ${C}_{11}(u),$ the ratio $\frac{{C}_{11}}{\hat{C}_1},$ is non-monotonic in $t$,\\ i.e.,
				$\frac{C_{11}(u)}{\prod_{i=1}^{n}u_{i}}$ is increasing for $0\leq t\leq F^{-1}(\min_{i} k_i),$ and  $\frac{C_{11}(u)}{C_{1}(u)}$ is decreasing for $t \geq F^{-1}(\max_{i} k_i),$ where $k_i=\frac{a_i-1}{a_i+b_i-1}$
\end{t1}
{\bf Proof.} Taking, $W=\frac{C_{11}(u)}{\prod_{i=1}^{n}u_{i}}$ we find $W=1+\alpha \Big(\prod_{i=1}^{n}u_{i}^{a_i-1}\Big)\Big( \prod_{i=1}^{n}(1-u_{i})^{b_i}\Big)=1+V,$ (say) respectively. We note that $\frac{dW}{dt}=\frac{dV}{dt}=Y\sum_{i=1}^{n}\Big(\frac{a_i-1-u_i(a_i+b_i-1)}{u_i(1-u_i)}\Big),$ and $\Big(\frac{a_i-1-u_i(a_i+b_i-1)}{u_i(1-u_i)}\Big),$ is positive if $0\leq u_i \in (0,k_i)$ and negative if $u_i \in (k_i,1),$ where $k_i=\frac{a_1-1}{a_i+b_i-1}$ for all $i=1,2,\ldots,n.$ Therefore, $\frac{dX}{dt}$ is positive  if  $0\leq t\leq F^{-1}(\min_{i} k_i),$ and negative if $t \geq F^{-1}(\max_{i} k_i).$ We  conclude that $\frac{C_{11}(u)}{\prod_{i=1}^{n}u_{i}}$ is increasing in $t$ for $0\leq t\leq F^{-1}(\min_{i} k_i),$  and decreasing  in $t$ for $t \geq F^{-1}(\max_{i} k_i).$$\hfill\Box$\begin{t1}
$\frac{C_{11}(u)}{\prod_{i=1}^{n}u_{i}}$ is non-monotonic in $t.$ Further, $\frac{\hat{C}_{11}}{\hat{C}_1},$ is non-monotonic in $t$,  \\
				$\frac{C_{11}(u)}{\prod_{i=1}^{n}u_{i}}\uparrow$ is for $0\leq t\leq F^{-1}(\min_{i} k_i),$ and  $\frac{C_{11}(u)}{C_{1}(u)}\downarrow$ for $t \geq F^{-1}(\max_{i} k_i),$ $k_i=\frac{a_1-1}{a_i+b_i-1}$  \end{t1}
    \subsection{   Linear Spearman (Bivariate copula)}
     Joe(1997), page 148 has defined the linear Spearman copulas as
     \begin{equation}\label{linear}
C(u_1,u_2)=\left\{ \begin{array}{lllll}
u_1u_2+\theta u_2 (1-u_1)& \textrm{if}&  u_2 \leq u_1, 0\leq \theta \leq 1,\nonumber\\
u_1u_2+\theta u_1 (1-u_2)& \textrm{if}&  u_2 > u_1, 0\leq \theta \leq 1,\nonumber\\
(1+\theta)u_1u_2 & \textrm{if}&  u_1+u_2 < 1, -1\leq \theta \leq 1,\nonumber\\
u_1u_2 +\theta (1-u_1)(1-u_2) & \textrm{if}&  u_1+u_2 \geq 1, -1\leq \theta \leq 0.\nonumber\\
\end{array} \right.
\end{equation}
It is easy to see the proof of following theorem.
    \begin{t1}
For Linear Spearman copula $C(u_1,u_2),$ the ratio $\frac{C(u_1,u_2)}{u_1u_2}$ is decreasing (increasing) in $t,$ according as $0\leq \theta \leq 1,$ and $-1\leq \theta \leq 0$.
\end{t1}
i.e., for Linear Spearman (Bivariate copula), we have
\begin{enumerate}
    \item [$(i)$]
 $T_{P}^{D}\leq_{rhr} T_{P}^{I}$ for $0\leq \theta \leq 1$;
 \item [$(ii)$] $T_{P}^{D}$ and  $T_{P}^{I}$ have proportional reversed hazard rate if $u_1+u_2<1$ and  $-1\leq \theta \leq 0$;
 \item[$(iii)$]$T_{P}^{D}\geq_{rhr} T_{P}^{I}$ if $u_1+u_2 \geq 1$ and $-1\leq \theta \leq 0.$
\end{enumerate}
In light of above theorem we state that $T_{P}^{D}\leq_{st} T_{P}^{I}$ and $T_{P}^{D}\leq_{mrl} T_{P}^{I}.$
A stronger result in terms of likelihood ratio order $lr$ is derived in the following theorem
\begin{t1}
If the reversed hazard rates $\mu_1$ and $\mu_2$ of $X_1$ and $X_2$ respectively are decreasing in $t$ then $\frac{f_{P}^{D}(t)}{f_{P}^{I}(t)}$ is decreasing in $t$ i.e., if $\mu_1$ and $\mu_2$ are decreasing in $t$ then $T_{P}^{D}\leq _{lr} T_{P}^{I}.$
\end{t1}
{\bf Proof.}
Here \begin{equation}
 f_{P}^{D}(t)
= (1-\theta) u_1u_2 \Big(\frac{u_2^{'}}{u_2}+\frac{u_1^{'}}{u_1}\Big)\nonumber
\end{equation}
where $'$
represents differentiation with respect to $t,$ and $f_{P}^{I}(t)=u_1u_2\Big(\frac{u_2^{'}}{u_2}+\frac{u_1^{'}}{u_1}\Big).$ \\Thus
\begin{eqnarray}
\frac{f_{P}^{D}(t)}{f_{P}^{I}(t)}&=&1-\theta+\frac{\theta}{u_1}\Big(\frac{u_2^{'}}{u_2}+\frac{u_1^{'}}{u_1}\Big), \nonumber\\
&=&1-\theta+\frac{\theta}{u_1}(\mu_1(t)+\mu_2(t))\nonumber
\end{eqnarray}
is decreasing in $t$ if the reversed hazard rates $\mu_1$ and $\mu_2$ of $X_1$ and $X_2$ respectively are decreasing in $t.$

    \begin{landscape}
\begin{table}
\caption{Error Analysis based on Copula functions (cf. Nadarajah et al. (2017)) with product copula, $C_1(u)=\prod_{i}^{n} u_{i}$}.\label{tab1}
  \scriptsize{  \begin{tabular}{|c|c|}
				\hline
				Copula/ Common Name by proposer & Stochastic ordering \\
				\hline
      $C_2(u)=\Big(\prod_{i}^{n} u_{i}\Big) \Big\{1+\alpha\prod_{i=1}^{n}\Big(1-u_{i}\Big)\Big\}$ & $\frac{C_2}{C_1}\downarrow (\uparrow)$ if $\alpha > (<) 0,$ i.e., $T_{P}^{D}\leq_{rhr} (\geq_{rhr}) T_{P}^{I}$ if $\alpha > (<) 0, \Rightarrow T_{P}^{D}\leq_{st} (\geq_{st})T_{P}^{I}$ if $\theta > (<) 0,$ \\
   Farlie-Gumbel-Morgenstern &$\frac{\hat{C}_2}{\hat{C}_1}\uparrow (\downarrow)$ if $\alpha > (<) 0,$ i.e., $T_{S}^{D}\geq_{hr} (\leq_{hr}) T_{S}^{I}$. So, $T_{S}^{D}\geq_{mrl} (\leq_{mrl}) T_{S}^{I}$,  $T_{S}^{D}\geq_{st} (\leq_{st}) T_{S}^{I}$ if $\alpha > (<) 0,$ \\\hline
   $C_3(u)=\Big(\prod_{i}^{n} u_{i}\Big) \Big\{1+\alpha \prod_{i=1}^{n}\Big(1-u_{i}^{1/r}\Big)\Big\}^{r}$ & $\frac{C_3}{C_1}\downarrow (\uparrow)$ if $\alpha > (<) 0,$ $T_{P}^{D}\leq_{rhr} T_{P}^{I}\Rightarrow T_{P}^{D}\leq_{st} T_{P}^{I}$\\
   Fischer and Kock &$\frac{\hat{C}_3}{\hat{C}_1}\downarrow (\uparrow)$ if $\alpha > (<) 0,$ $T_{S}^{D}\leq_{rhr} T_{S}^{I}\Rightarrow T_{S}^{D}\leq_{st} T_{S}^{I}$
    \\ \hline
   $C_4(u)=\Big\{\sum_{i=1}^{n}u_{i}^{-\alpha}-1\Big\}^{-\frac{1}{\alpha}}$ &$\frac{C_4}{C_1}\uparrow,$ i.e., $T_{P}^{D}\geq_{rhr} T_{P}^{I}\Rightarrow T_{P}^{D}\geq_{st} T_{P}^{I}$\\
				Clayton&$\frac{\hat{C}_{4}}{\hat{C}_1}\downarrow,$ i.e., $T_{S}^{D}\leq_{hr} T_{S}^{I},$ and hence $T_{S}^{D}\leq_{mrl} T_{S}^{I},$ $T_{S}^{D}\leq_{st}T_{S}^{I}$\\\hline
    $C_5(u)=\exp\Big[\Big\{-\sum_{i=1}^{n}\big(-\ln u_i\big)^{\alpha}\Big\}^{\frac{1}{\alpha}}\Big]$&$\frac{C_5}{C_1}\downarrow,$ i.e., $T_{P}^{D}\leq_{rhr} T_{P}^{I}\Rightarrow T_{P}^{D}\leq_{st} T_{P}^{I}$\\
	             Gumbel-Hougaard  &$\frac{\hat{C}_5}{\hat{C}_1}\uparrow,$ i.e., $T_{S}^{D}\geq_{hr} T_{S}^{I}$ giving $T_{S}^{D}\geq_{mrl} T_{S}^{I}$ and $T_{S}^{D}\geq_{st} T_{S}^{I}$\\\hline
$C_6(u)=\Big(\prod_{i=1}^{n}u_{i}\Big)\exp\Big(-\alpha \prod_{i=1}^{n}\ln u_i\Big)$&$\frac{C_6}{C_1}\downarrow,$ i.e., $T_{P}^{D}\leq_{rhr} T_{P}^{I}\Rightarrow T_{P}^{D}\leq_{st} T_{P}^{I}$\\
				Gumbel-Barnet  &\\\hline
   $C_7(u)=\frac{\prod_{i=1}^{n}u_{i}}{\Big\{1+\prod_{i=1}^{n}(1-u_i^{\alpha})\Big\}^{\frac{1}{\alpha}}}$ &$\frac{C_7}{C_1}\uparrow,$ i.e.,$T_{P}^{D}\geq_{rhr} T_{P}^{I}\Rightarrow T_{P}^{D}\geq_{st} T_{P}^{I}$\\
				Nelson Ten &$\frac{\hat{C}_7}{\hat{C}_1}\downarrow,$ i.e.,$T_{S}^{D}\leq_{hr} T_{S}^{I}$ giving $T_{S}^{D}\leq_{mrl} T_{S}^{I}$ and $T_{S}^{D}\leq_{st} T_{S}^{I}$\\\hline
$C_8(u)=\Big(\prod_{i=1}^{n}u_{i}^{\alpha_i}\Big)\min\Big\{u_1^{\alpha_1}, u_2^{\alpha_2}, \ldots, u_n^{\alpha_n} \Big\}$&$\frac{C_8}{C_1}\uparrow$ if $\alpha >1$, i.e., $T_{P}^{D}\geq_{rhr} T_{P}^{I}\Rightarrow T_{P}^{D}\geq_{st} T_{P}^{I}$\\
				Marshall and Olkin& $\frac{\hat{C}_8}{\hat{C}_1}5\downarrow$ if $\alpha >1$, i.e., $T_{P}^{D}\geq_{rhr} T_{P}^{I}\Rightarrow T_{P}^{D}\geq_{st} T_{P}^{I}$\\\hline
    $C_9(u)=(1-\alpha)\Big\{\prod_{i=1}^{n}\Big(\frac{1-\alpha}{u_i}+\alpha\Big)-\alpha\Big\}^{-1}$&$\frac{C_9}{C_1}\uparrow (\downarrow)$ if $\alpha \in [-1,0] (\alpha \in [0,1]),$ i.e., $T_{P}^{D}\geq_{rhr} (\leq_{rhr}) T_{P}^{I}\Rightarrow T_{P}^{D}\geq_{st} (\leq_{st})T_{P}^{I}$\\
				(AMH) Ali-Mikhail-Haq
    & $\frac{\hat{C}_9}{\hat{C}_1}\downarrow (\uparrow)$ if $\alpha \in [-1,0] (\alpha \in [0,1]),$ i.e., $T_{S}^{D}\leq_{hr} (\geq_{hr}) T_{S}^{I},$
    So, $T_{S}^{D}\leq_{mrl} (\geq_{mrl}) T_{S}^{I}$,  $T_{S}^{D}\leq_{st} (\geq_{st}) T_{S}^{I}$
\\\hline
    $C_{10}(u)=\Big[\Big\{\alpha\min_{i}u_i\Big\}^{m}+\Big\{(1-\alpha) \prod_{i=1}^{n}u_i\Big\}^{m}\Big]^{1/m}$ &$\frac{C_{10}}{C_1}\downarrow,$ i.e.,$T_{P}^{D}\leq_{rhr} T_{P}^{I}\Rightarrow T_{P}^{D}\leq_{st} T_{P}^{I}$\\
				Fischer and Hinzmann's copulas& $\frac{\hat{C}_{10}}{\hat{C}_1}\uparrow,$ i.e.,$T_{S}^{D}\geq_{hr} T_{S}^{I}$ giving $T_{S}^{D}\geq_{mrl} T_{S}^{I}$ and $T_{S}^{D}\geq_{st} T_{S}^{I}$\\\hline
    $C_{11}(u)=\prod_{i=1}^{n}u_{i}(t)+\alpha \Big(\prod_{i=1}^{n}u_{i}^{a_i}\Big)\Big( \prod_{i=1}^{n}(1-u_{i})^{b_i}\Big)$ &$\frac{\hat{C}_{11}}{\hat{C}_1},$ is non-monotonic in $t$,  \\
				Extension of Rodriguez-Lallena and Ubeda Flores& $\frac{C_{11}(u)}{\prod_{i=1}^{n}u_{i}}\uparrow$ is for $0\leq t\leq F^{-1}(\min_{i} k_i),$ and  $\frac{C_{11}(u)}{C_{1}(u)}\downarrow$ for $t \geq F^{-1}(\max_{i} k_i),$ $k_i=\frac{a_1-1}{a_i+b_i-1}$ \\\hline
    Linear Spearman (Bivariate copula)& for $0\leq \theta \leq 1,$ $\frac{\hat{C}_{10}}{\hat{C}_1}\downarrow,$ i.e.,
 $T_{P}^{D}\leq_{rhr} T_{P}^{I};$ \\
   & for $-1\leq \theta \leq 0,$ $\frac{\hat{C}_{10}}{\hat{C}_1}$ is constant ($T_{P}^{D}$ and  $T_{P}^{I}$ have PRHR) until $u_1+u_2<1,$\\& $\frac{C(u_1,u_2)}{u_1u_2}\uparrow$, i.e.,  $T_{P}^{D}\geq_{rhr} T_{P}^{I}$ whenever $u_1+u_2 \geq 1.$ \\\hline
   \end{tabular}}
	\end{table}
\end{landscape}
\section{Conclusion}
\begin{enumerate}
\item[(i)] Nanda et al. (2022) discussed the error essessment for Gumbel-I, Gumbel-II,  Gumbel-III distributions. We note that Gumbel-I, Gumbel-II, and Gumbel-III can be obtained from Gumbel Barnet, Farlie Gumbel  Morgenstern and Gumbel-Hougaard copula functions respectively by taking $\hat{u}_{i}(t)=\exp(-\lambda_{i}t),t \geq 0,$ and $i=1,2.$\\ ${}$\hspace{0.8cm}The stochastic orders depicted in Table \ref{tab1} related to the said copulas reestablish the same facts as stipulated by Nanda et al. (2002) about OA and UA in the reliability functions except in mean residual life function in Gumbel-I, survival function of Gumbel-II, for which a mismatch has been found and the same has been discussed in detail in the subsequent point.
\item[(ii)] We find a discrepancy in the results tabulated in Table 4 of Nanda et al. (2022), particularly with regard to Gumbel I and Gumbel II distributions. They mentioned that there is over assessment (OA) in hazard rate, i.e., $T^{D}_S \leq_{hr} T^{I}_S.$ As a result, it follows that $T_P \geq_{st} T_S,$ $T_P \geq_{mrl} T_S.$ So, we must have under-assessment (UA) in reliability and mean residual life function. So it is evident that the findings of  Nanda et al. (2022) with regard to error analysis in reliability function of Gumbel II are incorrect. On the other hand, if there is under-assessment in hazard rate for Gumbel-I distribution  then we should have over-assessment in mean residual life function. But the results in Table 4 of Nanda et al. for Gumbel-I do not depict the same as far as mean residual life function is concerned.
\end{enumerate}
 \begin{landscape}
\begin{table}
 \caption{Error assessment:  $C_1(u)=\prod_{i}^{n} u_{i}$}
 \label{tab1}
   \begin{tabular}{|c|c|c|}
				\hline
				Copula/ Common Name by proposer &   Parallel System&Series System\\
				\hline
      $C_2(u)=\Big(\prod_{i}^{n} u_{i}\Big) \Big\{1+\alpha\prod_{i=1}^{n}\Big(1-u_{i}\Big)\Big\}$   &$\alpha>0:~$rhr,sf (OA)&$\alpha>0:~$hr (OA),sf,mrl:UA\\
     farlie & rhr,sf (UA) &\\\hline
   $C_3(u)=\Big(\prod_{i}^{n} u_{i}\Big) \Big\{1+\alpha \prod_{i=1}^{n}\Big(1-u_{i}^{1/r}\Big)\Big\}^{r}$  &$\alpha<0:~$rhr,sf (UA)&$\alpha<0:~$hr,sf,mrl:UA\\
   Fischer and Kock &&
    \\ \hline
   $C_4(u)=\Big\{\sum_{i=1}^{n}u_{i}^{-\alpha}-1\Big\}^{-\frac{1}{\alpha}}$ &&\\
				Clayton&&\\\hline
    $C_5(u)=\exp\Big[\Big\{-\sum_{i=1}^{n}\big(-\ln u_i\big)^{\alpha}\Big\}^{\frac{1}{\alpha}}\Big]$&&\\
	             Gumbel-Hougaard  &&\\\hline
$C_6(u)=\Big(\prod_{i=1}^{n}u_{i}\Big)\exp\Big(-\alpha \prod_{i=1}^{n}\ln u_i\Big)$&&\\
				Gumbel-Barnet  & &\\\hline
   $C_7(u)=\frac{\prod_{i=1}^{n}u_{i}}{\Big\{1+\prod_{i=1}^{n}(1-u_i^{\alpha})\Big\}^{\frac{1}{\alpha}}}$ &&\\
				Nelson Ten &&\\\hline
$C_8(u)=\Big(\prod_{i=1}^{n}u_{i}^{\alpha_i}\Big)\min\Big\{u_1^{\alpha_1}, u_2^{\alpha_2}, \ldots, u_n^{\alpha_n} \Big\}$&&\\
				Marshall and Olkin& &\\\hline
    $C_9(u)=(1-\alpha)\Big\{\prod_{i=1}^{n}\Big(\frac{1-\alpha}{u_i}+\alpha\Big)-\alpha\Big\}^{-1}$&&\\
				(AMH) Ali-Mikhail-Haq
    & &
\\\hline
    $C_{10}(u)=\Big[\Big\{\alpha\min_{i}u_i\Big\}^{m}+\Big\{(1-\alpha) \prod_{i=1}^{n}u_i\Big\}^{m}\Big]^{1/m}$ &&\\
				Fischer and Hinzmann's copulas& &\\\hline
    $C_{11}(u)=\prod_{i=1}^{n}u_{i}(t)+\alpha \Big(\prod_{i=1}^{n}u_{i}^{a_i}(1-u_{i})^{b_i}\Big)$ &  &\\
				Extension of Rodriguez-Lallena &  &\\
    and Ubeda Flores &&\\\hline
    Linear Spearman (Bivariate copula)&  &\\
   & & \\\hline
   \end{tabular}
\end{table}
\end{landscape}

\section*{Acknowledgement}
Subarna Bhattacharjee  would like to thank Odisha State Higher Education Council for providing support to carry out the research project under OURIIP, Odisha, India (Grant No. 22-SF-MT-073).

\end{document}